\documentclass[12pt, draftcls, onecolumn]{IEEEtran}  

\IEEEoverridecommandlockouts                              
\usepackage{graphics,epsfig,epic,eepic,theorem,latexsym,amssymb,amsmath,picinpar}

\newcommand{\rr}{\mathop{{\rm I}\mskip-4.0mu{\rm R}}\nolimits}

\title{\LARGE \bf
Proofs of ''LQG Control For MIMO Systems Over Multiple TCP-like Erasure Channels''}
\author{E. Garone, B. Sinopoli, A. Goldsmith and A. Casavola
\thanks{E. Garone and A. Casavola are with the Dipartimento di Elettronica, Informatica e Sistemistica, Universit\`{a} degli Studi della Calabria, Via Pietro Bucci, Cubo 42-c, Rende (CS), 87036, Italy
        {\tt\small \{egarone,casavola\}@deis.unical.it}}%
\thanks{B. Sinopoli is with the Department of Electrical and Computer Engineering, Carnegie Mellon University,
        Pittsburgh, PA 15213, USA
        {\tt\small brunos@ece.cmu.edu}}%
\thanks{A. Goldsmith is with the Department of Electrical Engineering, Stanford University, Stanford, CA 15213, USA
        {\tt\small andrea@ee.stanford.edu}}%
}

\begin{document}
\maketitle
\thispagestyle{empty}
\pagestyle{empty}


Here we will provide the proofs of the results stated in the
Infinite Horizon LQG Control section of \cite{GaroneTAC}
 by focusing on the control law and the related MARE. The analysis of the observation case can be
achieved in a dual way and it is partially covered by \cite{Liu:04}.
\\
\\
Let us recall the control MARE:
\begin{equation}\label{MAREApp}
\begin{array}{rcl}
S_{k+1}&=&\Pi_c\left(S_{k},A,B,U,W,\bar{N}\right)\\
&=&\!A^T\!S_{k}A\!+\!W\!-\!A^T\!S_{k}B\bar{N}\! \left[\! {\sum\limits_{I\in 2^\Im } \left[\eta_I^2
\left( {{\rm N}_I \!\left( {U\!+\!B^T\!S_{k}B} \right){\rm N}_I } \right)\right]}
\right]^{-1}\!\!\!\!\!\!\!\bar{N} B^T\!S_{k}A,
\end{array}
\end{equation}
where $\eta_I\in \rr$ is
\[
\eta_I= \sqrt{\left( {\prod\limits_{i\in I}{{\bar{\nu} _i } } }
\right)\left( {\prod\limits_{i\notin I}{ {(1-\bar{\nu} _i)
} } } \right)},
\]
for every set $I \in \Im=\{1,...,m\}$ . Moreover, let us recall the auxiliary function
\begin{equation}\label{op1Bis2}
\phi \left( {K,X} \right)=\sum\limits_{I\in 2^\Im } \left[\eta_I^2
\left( {F_I XF_I ^T+V_I } \right)\right],
\end{equation}
where
\begin{eqnarray}\label{FIVI1}
 F_I &\stackrel{\Delta}{=}&A^T+K\left( {N_I B^T} \right), \\ \label{FIVI2}
 V_I &\stackrel{\Delta}{=}&W+KN_I UN_I ^TK^T.
\end{eqnarray}
The following results can be proved.
\\
\noindent {\bf Proposition (P1)} - If $K$ is chosen such that
\[K=\bar{K}_X=-A^TXB\bar{N} \left[ {\sum\limits_{I\in 2^\Im
} {\left[\eta_I^2 \left( {{\rm N}_I \left( {U+B^TXB} \right){\rm
N}_I } \right)\right]} } \right]^{-1},
\]
then
\[
\phi \left({\bar{K}_X,X} \right)= \Pi_c\left(X,A,B,U,W,\bar{N}\right).
\]
\noindent{\bf Proof} -
Consider the operator (\ref{op1Bis2}) together with (\ref{FIVI1})-(\ref{FIVI2})
\[
\phi \left( {K,X} \right)\!=\!\sum\limits_{I\in 2^\Im }\left[\eta_I^2 {
\left( {\left( {A^T\!\!+\!KN_I B^T\!} \right)X\left( {A\!+\!BN_I K^T}
\right)\!+\!W\!+\!KN_I UN_I K^T} \right)}
\right].
\]
Let us expand the product and group the terms of the latter equation
\[
\begin{array}{rcl}
\phi \!\left( {K,X} \right)\!&=&\!\sum\limits_{I\in 2^\Im }\bigg[  \eta_I^2
\left(\!A^T\!XBN_I K^T\!\!+\!KN_IB^T\!XA\!+KN_IB^T\!XBN_IK^T+A^T\!XA\!+W+\!\right.\\
& & \left. \!+\!KN_I UN_I K^T \right)\bigg]=\\
&=&\!\sum\limits_{I\in 2^\Im }\bigg[  \eta_I^2
\left(\!A^T\!XBN_I K^T\!\!+\!KN_IB^T\!XA\!+KN_I\left(B^T\!XB+U\right)N_IK^T\!\right)+\\
& & +\eta_I^2\left( \!A^T\!XA\!+W\! \right)\bigg]
\end{array}
\]
By splitting the summation and exploiting the fact that $\sum_{I\in 2^\Im }\eta_I^2=1,$ we obtain
\[
\begin{array}{rcl}
\phi \!\left( {K,X} \right)\!&=&\!\sum\limits_{I\in 2^\Im }\bigg[  \eta_I^2
\left(\!A^T\!XBN_I K^T\!\!+\!KN_IB^T\!XA\!+KN_I\left(B^T\!XB+U\right)N_IK^T\!\right)\bigg]+\\
& & +\!\left(\sum\limits_{I\in 2^\Im } \eta_I^2 \right)\left( \!A^T\!XA\!+W\! \right)=\\
&=&\!\sum\limits_{I\in 2^\Im }\!\bigg[  \eta_I^2\!
\left(\!A^T\!XBN_I K^T\!\!+\!KN_IB^T\!XA\!+\!KN_I\!\left(B^T\!XB+U\right)N_IK^T\!\right)\bigg]\!+\!A^T\!XA\!+W.\!
\end{array}
\]
Moreover, because $\sum_{I\in 2^\Im }\eta_I^2 N_I=\bar{N}$ we have
\[
\begin{array}{rcl}
\phi \!\left( {K,X} \right)\!&=&\!\sum\limits_{I\in 2^\Im }\!\bigg[  \eta_I^2\!
\left(\!\!KN_I\!\left(B^T\!XB+U\right)N_IK^T\!\right)\bigg]\!\!+A^T\!XB\bar{N} K^T\!\!+\!K\bar{N}B^T\!XA\!+A^T\!XA\!+\!W\\
&=&\!KL^{-1}(X)K^T\!\!+A^T\!XB\bar{N} K^T\!\!+\!K\bar{N}B^T\!XA\!+A^T\!XA\!+\!W.
\end{array}
\]
where $L(X)$ is a shorthand defined as follows
\[
L(X)=\left[ {\sum\limits_{I\in 2^\Im }\eta_I^2 {\left( {{\rm N}_I
\left( {U+B^TXB} \right){\rm N}_I } \right)} } \right]^{-1}.
\]
Finally, by substituting $K=\bar{K}_x=-A^TXB \bar{N} \left( L(X) \right)$ we obtain
\[
\begin{array}{l}
\phi\! \left( K,X \right)\!=\\
=\!A^T\!XB \bar{N} L(X)\bar{N}B^T\!SA\!-\!A^T\!XB \bar{N} L(X) \bar{N} B^T\!XA\!-\!A^T\!XB \bar{N} L(X) \bar{N} BXA+\!A^T\!XA\!+\!W\!=\\
=\!A^T\!XA\!+\!W\!-\!A^T\!XB \bar{N} L(X) \bar{N} B^T\!XA,
\end{array}
\]
which concludes the proof. \hfill $\Box$
\\
\\
%
%
\noindent {\bf Proposition (P2)} - Let us define
\begin{equation}\label{g_barN}
g_{\bar{N}}
\left( X \right)=\min _k \phi \left( {K,X} \right),
\end{equation}
then
\[
g_{\bar{N}} \left( X \right)=\Pi_c\left(X,A,B,U,W,\bar{N}\right).
\]
\noindent{\bf Proof} - By differentiate $\phi(K,X)$
with respect to $K$ we obtain:
\[
 \frac{\partial \phi \left( {K,X} \right)}{\partial
K}=2K(L(X))^{-1}+2A^TXB\bar{N} =0.
\]
The minimizer is then
\[
K=-A^TXB\bar{N} L(X)=\bar{K}_X,
\]
and because of Proposition 1, the statement follows. \hfill $\Box$
\\
\\
\noindent {\bf Proposition (P3)} - If $X\le Y$ then $g_{\bar{N}} \left( X \right)\le g_{\bar{N}} \left( Y
\right).$
\\
\normalsize \noindent{\bf Proof} - Being $\phi \left( {K,X} \right)$
affine in $X$, then\[g_{\bar{N}} \left( X
\right)=\mathop {\min }\limits_K \phi \left( {K_X ,X} \right)\le
\phi \left( {K_Y ,X} \right)\le \phi \left( {K_Y ,Y}
\right)=g_{\bar{N} } \left( Y \right).\] \hfill $\Box$
\\
%
\noindent {\bf Lemma (L1)} - Let us define the following operator:
\[
{\cal{L}} \left( Y \right)=\sum\limits_{I\in 2^\Im } {\eta_I^2\ \left( {F_I YF_I ^T} \right)},
\]
If there exists a positive matrix $\bar{Y} >0$ such that $\bar{Y} >{\cal{L}}\left({\bar{Y}} \right)$
then:
\\
a) $\forall W\ge 0,\,\,\,\,\,\mathop {\lim }\limits_{k\to \infty }
{\cal{L}}^k\left( W \right)=0,$
\\
b) Let  $U\ge 0$ and let $Y_{k+1} ={\cal{L}}\left( {Y_k } \right)+U$
initialized at Y$_{0}$: then the sequence Y$_{k}$ is bounded.
\\
\normalsize \noindent{\bf Proof} -
\\
a) Let us choose two scalars $r\in [0,1)$ and $m \ge 0$ such that
\begin{itemize}
\item ${\cal{L}}\left( {\bar{Y} } \right)<r\bar{Y},$
\item $W\le m\bar{Y}.$
\end{itemize}
Being the operators ${\cal{L}}$ linear, crescent monotone, i.e. if $X>Y$ then ${\cal{L}}\left( X \right)>{\cal{L}}(Y)$,
and because $W \geq 0$ implies ${\cal{L}}(W) \geq 0$, the following inequality results
\[
0 \leq {\cal{L}}^k\left( W \right)\le m{\cal{L}}^k\left( {\bar{Y} }
\right)<mr^k\bar{Y}.
\]
Then, for $k\to \infty$, we have ${\cal{L}}^k\left( W \right)\to 0.$
\\
b) Let introduce two further scalars $m_U \ge 0$, $m_{Y_0} \ge 0$. By following the same lines of a) we obtain the following inequality
\[
Y_k ={\cal{L}}^k\left( {Y_0 } \right)+\sum\limits_{t=0}^{k-1}
{{\cal{L}}^t\left( U \right)} < \left( {m_{Y_0 }
r^k+\sum\limits_{t=0}^{k-1} {m_U r^t} } \right)\bar{Y} =
\left( {m_{Y_0 } r^k+m_U\frac{1-r^k }{1-r}} \right)\bar{Y},
\]
where the properties of the geometric series have been exploited. \hfill $\Box$
\\
\\
\noindent {\bf Lemma (L2)} - Let us suppose there exits a pair $(\bar{K},\bar{S})$ of matrices such that:
\[
\begin{array}{rcl}
\bar{S} &\ge& 0,\\
\bar{S} &>& \phi \left( {\bar{K},\bar{S} } \right),
\end{array}
\]
then, $\forall S_0 $ the sequence $S_k =g_{\bar{N}}^k \left( {S_0 }
\right)$ is bounded.
\\
\noindent {\bf Proof} - By using the operator
\[
{\cal{L}}\left( Y \right)=\sum\limits_{I\in 2^\Im } \eta_I^2
\left( {F_{I} YF_{I} ^T} \right),
\]
where
\[
\begin{array}{l}
 F_{I} \stackrel{\Delta}{=}A^T+\bar{K}\left( {N_I B^T} \right), \\
 \end{array}
\]
we can write
\[
\bar{S} >\phi \left( {\bar{K},\bar{S} } \right)={\cal{L}}\left(
{\bar{S} } \right)+W+\sum\limits_{I\in 2^\Im } {\eta_I^2 \left( {\bar{K}N_I UN_I ^T\bar{K}^T} \right)\ge } {\cal{L}}\left(
{\bar{S} } \right).
\]
By exploiting the definition of $g(S_k)$
\[
S_{k+1} =g\left( {S_k } \right)\le \phi \left( {K,S_k }
\right)={\cal{L}}\left( {S_k } \right)+W+\sum\limits_{I\in 2^\Im }
{\eta_I^2 \left( {\bar{K}N_I UN_I ^T\bar{K}^T} \right)}
={\cal{L}}\left( {S_k} \right)+U,
\]
where
\[
U=W+\sum\limits_{I\in 2^\Im } {\eta_I^2\left( {\bar{K}N_I UN_I
^T\bar{K}^T} \right)} \ge 0.
\]
Then by resorting to Lemma L1 we can conclude that {\{}S$_{k}${\}} is bounded. \hfill $\Box$
\\
\\
\noindent {\bf Lemma (L3)} - Let $X_{t+1} =h\left( {X_t }
\right),Y_{t+1} =h\left( {X_t } \right)$, if $h$ is a monotone function, then:
\[
\begin{array}{l}
 X_1 \le X_0 \Rightarrow X_{t+1} \le X_t \,\,\forall t, \\
 X_1 \ge X_0 \Rightarrow X_{t+1} \ge X_t \,\,\forall t, \\
 Y_0 \le X_0 \Rightarrow Y_t \le X_t \,\,\forall t. \\
 \end{array}
\]
\noindent {\bf Proof} -
Consider the first condition. It is true t=0. Then, by induction:
\[
X_{t+2} = h\left( {X_{t+1} } \right)\le h\left( {X_t } \right)=X_t
\]
The other two cases are similar. \hfill $\Box$
\\
\noindent {\bf Lemma (L4)} -  Consider the operator (\ref{op1Bis2}) together with (\ref{FIVI1})-(\ref{FIVI2})
If it exists a couple of matrices $(\tilde{K} ,\tilde{S})$ such that
\[
\begin{array}{rcl}
\tilde{S} & > & 0, \\
\tilde{S} & > & \phi \left(
{\tilde{K} ,\tilde{S} } \right),
\end{array}
\]
then:
\begin{enumerate}
\item for each initial condition ${\bar{S}_0 } \ge 0$, the MARE (\ref{MAREApp}) converges and moreover
\[\mathop {\lim }\limits_{t\to \infty } S_k =\mathop {\lim }\limits_{k\to \infty } g_{\bar{N}}^k\left( {S_0 } \right)=\bar{S} \]
does not depend on the initial condition,
\item $\bar{S} $ is the unique positive definite fixed point of the MARE (\ref{MAREApp})
\end{enumerate}
\noindent {\bf Proof} -
\\
1) First let us prove that the MARE initialized at $S_{0}=Q_0=0$
converges. Let $Q_k =g_{\bar{N} }^k \left( 0 \right)$. Because Q$_{0}\le $ Q$_{1}$, from Proposition 3 it follows
\[
Q_1 =g_{\bar{N} } \left( {Q_0 } \right)\le g_{\bar{N} }
\left( {Q_1 } \right)=Q_2.
\]
By induction we have a monotonic nondecreasing sequence of matrix.
Because of Lemma L2 this sequence is be bounded as follows
\[
0\le Q_1 \le Q_2 \le ...\le M_{Q_0 }.
\]
Then, the sequence converges to one of the positive semidefinite fixed point of the MARE.
\[\mathop {\lim }\limits_{k\to \infty } Q_k =\bar{S},\]
such that $\bar{S} =g_{\bar{N} } \left( {\bar{S} } \right).$
\\
The next step is to show that the sequence converges to the same
point $\forall S_0=R_0 \ge \bar{S} $.
\\
By resorting to Lemma L2 notation, we can observe that
\[
\bar{S} =g_{\bar{N}}\left(\bar{S}\right)={\cal L}\left( {\bar{S} }
\right)+W+\sum\limits_{I\in 2^\Im } {\eta_I^2\left(
\bar{K}N_I UN_I ^T\bar{K}^T \right)}
>{\cal L}\left( {\bar{S}} \right)
\]
The latter implies, as a consequence of Lemma L1, that $\mathop {\lim }\limits_{k\to \infty } {\cal
L}^k\left( X \right)=0, \forall X \ge 0$.
\\
Let us assume $S_0=R_0 \ge \bar{S}$. Then
\[{\cal
L}\left( {R_0 } \right)=R_1 \ge {\cal
L}\left( {\bar{S}}
\right)=\bar{S}.
\]
By induction it yields $R_k \ge \bar{S} ,\,\,\forall k$.
\\
Notice that:
\[
\begin{array}{rcl}
0&\le& R_{k+1} -\bar{S} =g_{\bar{N}} \left( {R_k }
\right)-g_{\bar{N} } \left( {\bar{S} } \right)=\phi \left(
{\bar{K}_{R_k} ,R_k } \right)-\phi \left(
{\bar{K},\bar{S} } \right)\le \\
&\le& \phi \left(
{\bar{K} ,R_k } \right)-\phi \left(
{\bar{K} ,\bar{S}} \right)
=\sum\limits_{I\in 2^\Im } {\eta_I^2 \left(
{F_{I} \left( {R_k -S}
\right)F_{I}}
^T+V_{I} -V_{I}\right)={\cal L}\left( {R_k
-\bar{S}} \right)}.
\end{array}
\]
Then, finally
\[
0\le R_k-\bar{S} \le \mathop {\lim }\limits_{k\to \infty } {\cal
L}\left( {R_k -\bar{S}} \right)=
0.
\]
The last thing to show is that the Riccati euqtion converges to $\bar{S},\, \forall \bar{S}_0 \ge 0.$ Let us define
\[
\begin{array}{l}
 Q_0 =0, \\
 R_0 =S_0 +\bar{S}. \\
 \end{array}
\]
Then, $Q_0 \le S_0 \le R_0 $. By resorting to Lemma L3 we obtain:
\[
Q_k \le S_k \le R_k \,\,\,\,\forall k.
\]
Then, finally
\[
\bar{S} =\mathop {\lim }\limits_{k\to \infty } Q_k \le \mathop
{\lim }\limits_{k\to \infty } S_k \le \mathop {\lim }\limits_{k\to
\infty } R_k =\bar{S}.
\]
2) Let consider a certain $S'=g_{\bar{N}} \left(
S'\right)$. If the Riccati equation is initialized at $S_0 =S',$ a constant sequence will results.
Because each Riccati equation converges to $\bar{S} $, then $S'=\bar{S}$. \hfill
$\Box$
\\
\noindent {\bf Theorem (LMI)} - The following statements are
equivalent:
\\
1) $\exists \left( {\bar{K} ,\bar{S} } \right):\bar{S}
> 0\,\,,\,\,\,\bar{S} > \phi \left( {\bar{K} ,\bar{S} } \right)$
\\
2) $\exists Z, Y>0$ such that:
\[
\left[ {{\begin{array}{*{20}c}
 Y \hfill & Y \hfill & {\eta_\emptyset \left( {YA^T\!\!+\!ZN_\emptyset B}
\right)} \hfill & {\eta_\emptyset  ZN_\emptyset U^{1/2}} \hfill
& {...} \hfill & {\eta_\Im \left( {YA^T\!\!+\!ZN_\Im B} \right)}
\hfill & {\eta_\Im ZN_\Im U^{1/2}} \hfill \\
 Y \hfill & {W^{-1}} \hfill & 0 \hfill & 0 \hfill & {...} \hfill & 0 \hfill
& 0 \hfill \\
 {...} \hfill & 0 \hfill & {Y} \hfill & 0 \hfill & \hfill & {...}
\hfill & {...} \hfill \\
 \ast \hfill & {...} \hfill & 0 \hfill & I \hfill & \hfill & \hfill & \hfill
\\
 \ast \hfill & \hfill & {...} \hfill & {...} \hfill & \hfill & {...} \hfill
& \hfill \\
 \ast \hfill & \hfill & \hfill & \hfill & \hfill & {Y} \hfill & {...}
\hfill \\
 0 \hfill & {...} \hfill & \hfill & \hfill & {...} \hfill & 0 \hfill & I
\hfill \\
\end{array} }} \right]>0
\]
where the sequence $\emptyset $,{\ldots},$\Im $ represents all possible
sets belonging to 2$^{\Im }$
\\
\noindent {\bf Proof} -
Let us consider conditions 1) in the statement:
\[
\begin{array}{l}
S>0, \\
S>\phi \left( {K,S} \right)=\sum\limits_{I\in 2^\Im } {\eta_I^2 \left( {F_I XF_I ^T+V_I } \right)}, \\
 \end{array}
\]
where $F_I$ and $V_I$ are introduced in (\ref{FIVI1})-(\ref{FIVI2}).
By expanding the second condition, we obtain
\[
 S>\sum\limits_{I\in 2^\Im } {\eta_I^2 \left( {\left( {A^T+KN_I
B^T} \right)S\left(
{A^T+KN_I B^T} \right)^T+W+KN_I UN_I ^TK^T} \right)}
\]
which is equivalent to
\[
S-W-\sum\limits_{I\in 2^\Im } {\eta_I^2 \left( {\left( {A^T+KN_I
B} \right)S\left(
{A^TKN_I B^T} \right)^T+KN_I UN_I ^TK^T} \right)} >0 \\
\]
By iteratively applying the Schur's complements we obtain the following Matrix Inequality
\small
\[
\left[ {{\begin{array}{*{20}c}
 S \hfill & I \hfill & {\eta_\emptyset \left( {A^T\!\!+\!KN_\emptyset B}
\right)} \hfill & {\eta_\emptyset  KN_\emptyset U^{1/2}} \hfill
& {...} \hfill & {\eta_{\Im} \left( {A^T\!\!+\!KN_{\Im} B} \right)}
\hfill & {\eta_{\Im} KN_{\Im} U^{1/2}} \hfill \\
 \ast \hfill & {W^{-1}} \hfill & 0 \hfill & 0 \hfill & {...} \hfill & 0
\hfill & 0 \hfill \\
 {...} \hfill & 0 \hfill & {S^{-1}} \hfill & 0 \hfill & \hfill & {...}
\hfill & {...} \hfill \\
 \ast \hfill & {...} \hfill & 0 \hfill & I \hfill & \hfill & \hfill & \hfill
\\
 \ast \hfill & \hfill & {...} \hfill & {...} \hfill & \hfill & {...} \hfill
& \hfill \\
 \ast \hfill & \hfill & \hfill & \hfill & \hfill & {S^{-1}} \hfill & {...}
\hfill \\
 0 \hfill & {...} \hfill & \hfill & \hfill & {...} \hfill & 0 \hfill & I
\hfill \\
\end{array} }} \right]>0,
\] \normalsize
that, by means of the congruence transformation $diag{\{}S^{-1},I,I,{\ldots},I{\}}$, simplifies into
\small
\[
\left[\!\! {{\begin{array}{*{20}c}
 {S^{-1}}\!\! \hfill & {S^{-1}}\!\! \hfill & {{\eta_\emptyset }\!\! \left(
{S^{-1}\!A^T\!\!+\!S^{-1}\!KN_\emptyset B} \right)\!\!} \hfill & {
{\eta_\emptyset } S^{-1}\!KN_\emptyset U^{1/2}\!} \hfill & {...} \hfill &
{{\eta_\Im } \!\left( {S^{-1}\!A^T\!\!+\!S^{-1}\!KN_\Im B} \right)} \hfill &
{{\eta_\Im } S^{-1}\!KN_\Im
U^{1/2}\!} \hfill \\
 {S^{-1}\!\!} \hfill & {W^{-1}\!\!} \hfill & 0 \hfill & 0 \hfill & {...} \hfill & 0
\hfill & 0 \hfill \\
 {...} \hfill & 0 \hfill & {S^{-1}\!} \hfill & 0 \hfill & \hfill & {...}
\hfill & {...} \hfill \\
 \ast \hfill & {...} \hfill & 0 \hfill & I \hfill & \hfill & \hfill & \hfill
\\
 \ast \hfill & \hfill & {...} \hfill & {...} \hfill & \hfill & {...} \hfill
& \hfill \\
 \ast \hfill & \hfill & \hfill & \hfill & \hfill & {S^{-1}\!} \hfill & {...}
\hfill \\
 0 \hfill & {...} \hfill & \hfill & \hfill & {...} \hfill & 0 \hfill & I
\hfill \\
\end{array} }} \!\!\right]>0.
\]
\normalsize By substituting $Y=S^{-1}$,$Z=S^{-1}K$ and by noticing that
\[
S>0 \Leftrightarrow S^{-1}=Y >0,
\]
the proof is completed.
\hfill $\Box$

%
%


\end{document}